\documentclass[12pt]{amsart}
\usepackage{graphicx, amsfonts, amsthm, amsxtra, amssymb, verbatim}
\usepackage[mathscr]{euscript}

\newtheorem *{thm}{Theorem}
\newtheorem{EX}{Example}
\newtheorem{coro}{Corollary}

\newtheorem{lem}{Lemma}
\newtheorem{PROP}{Proposition}

\begin{document}

%\begin{center}
%{\Large \bf Embedded Curves and Foliations}
%\footnote{ %Math. classification:  \\
%Keywords: %}}
%\\
%\end{center}

%\vspace{0cm}
%\begin{center}

%{\large \sc Hossein Movasati and  Paulo Sad}  %\sc%
%Instituto de Matem\'atica Pura e Aplicada, IMPA, \\
%Estrada Dona Castorina, 110,\\
%22460-320, Rio de Janeiro, RJ, Brazil, \\
%E-mail:
%{\tt hossein@impa.br, sad@impa.br} %\end{center}

\title{Embedded Curves and Foliations}
\date{}
\author{Hossein Movasati and Paulo Sad}
%\author{Paulo Sad}

%\vspace{1cm}
\begin{abstract}
We prove the existence of regular foliations with a prescribed
tangency divisor in neighborhoods of negatively embedded holomorphic
curves; this is related to a linearization theorem due to Grauert.
We give also examples of neighborhoods which can not be linearized.
\end{abstract}

\maketitle

\vspace{1cm}

We consider in this paper the problem of finding regular holomorphic
foliations in neighborhoods of  smooth, compact, holomorphic  curves
embedded in  complex surfaces. More precisely, we fix a positive divisor of a curve
and ask whether there exists a holomorphic 
foliation whose divisor of tangencies with the curve is exactly that
divisor. Let us state our main result:
 
\begin {thm}
\label{main1} Let $C\hookrightarrow S$ be an embedding of the curve
$C$ into the surface $S$ such that $C\cdot C<0$.
\begin {itemize}
 
\item if $C\cdot C< 4-4g$, there exists a regular foliation defined in a neighborhood of
$C$ and transverse to $C$.
\item let a divisor $D= \sum_{k=1}^{l} n_kp_k$ be given in $C$, with $n_k \in {\mathbb N}_{>0}$. If 
$$
C\cdot C < 4-4g + \sum_{k=1}^{l}(n_k-1)
$$
there exists a regular foliation $\mathcal F$ defined in a neighborhood of $C$
which is transversal to $C$ except at the points $p_1,\cdots, p_l
\in C$, where ${\rm tang}_{p_k}(\mathcal F, C)= n_k$   for every
$k=1,\dots,l$.
\end {itemize}  
\end{thm}
 
\medskip

In the statement $C\cdot C$ stands for the self-intersection number of $C$ in $S$ or,
equivalently, the Chern class of the normal bundle of $C$ in $S$; the number $g$ is the genus of $C$. Each number 
${\rm tang}_{p_k}(\mathcal F, C)$ is the order of tangency at $p_k\in C$ between $C$ and the leaf of $\mathcal F$ that passes through the point $p_k$.

\bigskip
Our method to prove this theorem consists in i) find a holomorphic
line field defined {\it along} the curve $C$ with the prescribed set
of tangencies and the prescribed order of tangencies; for this
purpose we have no need to assume that the curve is negatively
embedded; ii) extend the line field to a neighborhood of the curve;
here we must work under the hypothesis  $C\cdot C<0$ in order to
assure the annihilation of some cohomology groups.

\medskip
We should mention that our primary motivation stems from a
Linearization Theorem due to Grauert(\cite{Gr}): a
curve possesses a neighborhood diffeomorphic to a neighborhood  
of the zero section of its normal bundle if the
embedding  is sufficiently negative ($C\cdot C< {\rm min}\{0, 4-4g\}$).  A proof can be done in two steps.
We start by  guaranteeing  the existence of a foliation transverse to the curve; this is the first case  of  our Theorem. Once this is acomplished the rest of the proof goes as in (\cite {CMS})by 
finding another holomorphic foliation in a neighborhood $V$ of $C$
which has $C$ as a leaf; this foliation and the transverse
one are used as a kind of coordinate system for $V$ when   
the desired diffeomorphism is constructed.

\medskip
In this paper we discuss also how to produce examples of embeddings such that
there are no foliations with a given divisor of tangencies when the
negativity condition is violated. In particular, examples where
linearization is not possible are presented. All these examples
depend on properties of line fields defined along the curve.

\section {\bf Line Fields and Embeddings}

Let us consider an embedding $C\hookrightarrow S$ of the compact,
smooth, holomorphic curve $C$ into the surface $S$. In this Section
we study existence of line fields defined along $C$; {\bf we do not assume
$C\cdot C<0$}.   Existence of a
line field  with a given divisor of tangencies is always granted
when the degree of
the divisor is sufficiently bigger then $C\cdot C$. On the other hand, uniqueness (but
perhaps not the existence) follows when this degree is not too big,
and we will see later how this leads to the construction of
interesting examples.

\bigskip

A holomorphic subbundle  
$Y\hookrightarrow TS|C$ is a {\it holomorphic line field along $C$}.
Equivalently we may say that a line field is a section of the
$\mathbb P^1$-bundle $\mathbb P(TS|C)$ over $C$. $Y$ has   a {\it
tangency} with $C$ at the point $p\in C$ when the morphism of line
bundles $Y\rightarrow NC= \dfrac{TS|C}{TC}$ has a zero at $p$; 
the order of the zero is the {\it
order of tangency} between $Y$ and $C$. We write the set of
tangencies as an effective divisor $D=\sum_{k=1}^l n_kD_k$ of $C$; the point
$p_k$ is a 
point of tangency of order $n_k$.

\smallskip
In order to motivate the next Proposition, let us remark that when
$Y$ is a line field along $C$ whose divisor of tangencies with $C$
is $D$ then $Y \simeq \mathcal O(-D)\otimes NC$ as line bundles. In
fact, the morphism $ Y\rightarrow NC$ seen as a section of
$H^0(C,Y^*\otimes NC)$ has $D$ as its divisor of zeroes; therefore
$Y^*\otimes NC\simeq \mathcal O(D)$. This allows us to confound a
line field along $C$ having $D$ as divisor of tangencies  with an
injective morphism $\mathcal O(-D)\otimes NC \rightarrow TS|C$. We will
from now use $c(NC)$ to denote the Chern class of the normal bundle $NC$
of $C$ in $S$; it is well know that $C\cdot C= c(NC)$.

\begin{PROP}
Let $D$ be an effective divisor of $C$, and assume
$$
 C\cdot C <4-4g+\sum_{k=1}^l (n_k-1)
$$
There exists an  injective bundle morphism $Y: \mathcal O(-D)\otimes
NC \rightarrow TS|_{C}$ which has $D$ as divisor of tangencies with
$C$.
\end{PROP}
\proof Let us use $L:= \mathcal O(-D)\otimes NC$ for simplicity.
Firstly we construct $Y$ locally, that is, in the restriction of  the line bundle
$L$  to a small open set  $U\subset C$ where it is isomorphic to $U\times {\mathbb C}$.
More precisely, we have i) an open subset
$U^{\prime}$ of $S$ with coordinates $(z_1,z_2)\in \mathbb C \times \mathbb C$ such
that $U= U^{\prime}\cap C$ is $\{z_2=0\}$; ii) a holomorphic
function $f$ of $U$ such $D|_U= \{f=0\}$ and iii) trivialization
coordinates $(z_1,t)$ for $L|_U$; then  we may define
$$
Y(z_1)(t)=(z_1,t,f(z_1)t).
$$
This can be done in each set of  an open covering $\{U_i\}_{i\in I}$ of $C$, so we get morphisms  $Y_i: L|{U_i}\rightarrow TS|{U_i}$ with the desired
property; we assume that the support of
each  $D|_{U_i}$ consists of a point at most and that there are
no points of tangency in the intersections $U_i\cap U_j$ when $i\neq j$.  Let
$\tilde Y_i$ denote the composition $L|{U_i}\rightarrow
TS|{U_i}\rightarrow NC|{U_i}$. As $ \tilde Y_i = a_{ij}\tilde Y_j $,
where  $\{a_{ij}\}\in H^1(C, {\mathcal O}^*(C))$ defines a line
bundle $J$, $ \{\tilde Y_i\}$ is a section of $J\otimes
Hom(L,NC)\simeq J\otimes L^*\otimes NC$ having $D$ as divisor of
zeroes, so that $J\otimes L^*\otimes NC\simeq \mathcal O(D)$.
Consequently $J$ is the trivial line bundle and  we may suppose
$a_{ij}=1$, or $\tilde Y_i=\tilde Y_j$.

Now we have that 
$$ 
\{Y_{ij}\}:= \{Y_i-Y_j\}\in H^1(C,
Hom(L,TC))\simeq H^1(C,L^*\otimes TC);
$$

\noindent Let $\tilde D = \sum_{k=1}^l p_k$ and $s=\{s_i\}\in H^0(C,\mathcal O(\tilde D))$ whose
divisor of zeroes is $\tilde D$. Therefore
$$
(Y_i-Y_j)\otimes s^{-1}\in H^1(C,\mathcal O(-{\tilde D})\otimes L^*\otimes TC)
$$
\noindent and by Serre's duality 
$$
H^1(C,\mathcal O(-{\tilde D})\otimes L^*\otimes TC)\simeq H^0(C,KC^2\otimes \mathcal O(\tilde D)\otimes \mathcal O(-D)\otimes NC)
$$
\noindent($KC$ stands for the canonical bundle of $C$). By hypothesis the Chern class of the line bundle $KC^2\otimes \mathcal O(\tilde D)\otimes \mathcal O(-D)\otimes NC$ is negative; we conclude that $(Y_i-Y_j)\otimes s^{-1}= X_i-X_j$ for $X_i\in H^0(U_i,\mathcal O(-{\tilde D})\otimes L^*\otimes TC)$, and therefore $Y_i-Y_j=(X_i-X_j)\otimes s= s_iX_i-s_jX_j$. We define $Y:= Y_i-s_iX_i$ in each $U_i$. Clearly $Y$ is injective outside the support of $\tilde D$;
at each $p_i$, it is equal to $Y_i$, so it is also injective. As for the order of tangency at a point $p_i$, it coincides with the order of tangency of $Y_i$, which is $n_i$ by construction. \qed

\bigskip Consequently, there exists always a holomorphic line field along any curve
  if we admit a  number of
tangencies sufficiently big. We see also that there exists always a
holomorphic line field with any number of tangencies if $C\cdot C< 4-4g$.

\bigskip
\begin{PROP}
 \label{25mar11}
Assume that $C\cdot C <4-4g$.
There exists an  injective bundle morphism $Y: NC \rightarrow TS|_{C}$ which has no tangencies with $C$.
\proof We just have to repeat the  arguments   applied above  without
the presence of tangencies. We see that the condition in the statement
implies that $Y_i-Y_j = X_i-X_j$ for $X_i \in H^0(U_i,NC^*\otimes TC)$. \qed

\end{PROP}

In the next section we will analyse how to extend this
holomorphic line field to a neighborhood of the curve. For the
moment, let us state a general result concerning uniqueness.

\begin{PROP} 
\label{25mar2011} 
\label {atmostone}
Let D be an effective divisor of $C$ and assume
$$
%H^0(C, \mathcal O(D)\otimes NC^*\otimes TC)=0
c(NC)>2-2g+\sum n_i
$$
There exists at most one   line field    along $C$ having $D$ as
divisor of tangencies.
\end{PROP}
\proof Let us consider  two such line fields $Y_1$ and $Y_2$ as
bundle morphisms from $\mathcal O(-D)\otimes NC$ into $TS|_C$. The
induced morphisms $\tilde Y_i:\mathcal O(-D)\otimes NC\rightarrow NC$
seen as sections of $\mathcal O(D)\otimes NC^*\otimes NC= \mathcal
O(D)$ have the same divisor $D$ of zeroes, so that $\tilde Y_1=a
\tilde Y_2$ for some $a\in \mathbb C^*$. It follows that $Y_1-aY_2$
is a bundle morphism from $\mathcal O(-D)\otimes NC$ to $TC$; the
hypothesis tells us that $\mathcal O(D)\otimes
NC^*\otimes TC$ is a negative line bundle and so $Y_1-aY_2=0$. \qed

\medskip

%%%%%%%%%%%%%%%%%%%%%%%%%%%%%%%%%%%%%%%%%%%%%%%%%%%%%5
\smallskip
\section {\bf Neighborhoods of Negatively Embedded Curves}
Before proving the Theorem stated in the Introduction, we collect
some properties due to Grauer that are verified in the case of a negatively
embedded curve $C\hookrightarrow S$ (\cite{CM},\cite {Gr}).

\begin{itemize}

\item $C$ has a fundamental system of strictly pseudoconvex neighborhoods in $S$.

\item if $\mathcal G$ is a coherent sheaf defined in one of these neighborhoods, say
$V$, and $\mathcal I_C$ is the ideal sheaf of $C$ in $V$ then $$
\exists k>0 \,\,\,{\rm such\,\,\, that}\,\,\, H^i(V, {\mathcal
I}_C^k\cdot {\mathcal G})=0,\,\,\ i=1,2.
$$
\end {itemize}
\bigskip

\begin {lem} We have $H^2(V, {\mathcal I}_C \cdot {\mathcal G})=0$. Moreover if $$
H^0(C,KC\otimes NC^{\nu}\otimes {\mathcal G}^{*}|C)=0 $$ for all
$\nu \geq 1$ then $H^1(V, {\mathcal I}_C\cdot {\mathcal G})=0$.
\end{lem}

\proof From $H^i(V,{\mathcal I}_C^{\nu}/{\mathcal I}_C^{\nu +1}\cdot
\mathcal G)\simeq H^i(C, (NC^*)^{\nu}\otimes {\mathcal G}|C)$ we
get immediately $H^2(V,{\mathcal I}_C^{\nu}/{\mathcal I}_C^{\nu
+1}\cdot{\mathcal G})=0$. As  $$ H^1(C, (NC^*)^{\nu}\otimes
{\mathcal G}|C) \simeq H^0(C,KC\otimes NC^{\nu}\otimes {\mathcal
G}^{*}|C)
$$ (by Serre's duality)  we get $H^1(V,{\mathcal I}_C^{\nu}/{\mathcal I}_C^{\nu
+1}\cdot{\mathcal G})=0$ as well.

Let us consider the short exact sequence

$$
0 \rightarrow {\mathcal I}_C^{\nu +1}\cdot{\mathcal G} \rightarrow
{\mathcal I}_C^{\nu}\cdot{\mathcal G} \rightarrow {\mathcal
I}_C^{\nu}/{\mathcal I}_C^{\nu +1}\cdot {\mathcal G}\rightarrow 0
$$

\noindent which leads to

$$
\cdots \rightarrow H^i(V,{\mathcal I}_C^{\nu + 1}\cdot{\mathcal
G})\rightarrow H^i(V, {\mathcal I}_C^{\nu}\cdot{\mathcal
G})\rightarrow H^i(V,{\mathcal I}_C^{\nu}/{\mathcal I}_C^{\nu
+1}\cdot{\mathcal G})\rightarrow \cdots 
$$

\bigskip
\noindent Therefore the maps
$H^i(V,{\mathcal I}_C^{\nu + 1}\cdot{\mathcal G})\rightarrow H^i(V,
{\mathcal I}_C^{\nu}\cdot{\mathcal G}),\ \ i=1,2$, are  always
surjective. Consequently  $H^i(V,{\mathcal I}_C^k \cdot \mathcal G)
= 0$ for some $k>0$ implies $H^i(V, {\mathcal I}_C \cdot \mathcal
G)=0$, $i=1,2$. \qed

\bigskip
The next Lemma allows us to extend any line bundle over $C$ to a
line bundle over $V$. Of course there are certain line bundles which
are extendible regardless of the negativity of the embedding
$C\hookrightarrow V$. For example, $KC = KV|C \otimes NC=
KV|C\otimes [C]|C$, so that $KC$ always has an extension to $V$.
Below in our Theorem we find this situation when no tangencies are
present.

\begin {lem} The restriction $H^1(V,{\mathcal O}_V^{*})\rightarrow H^1(C,{\mathcal
O}_C^{*})$ is surjective.
\end {lem}
\proof Let $J$ be the subsheaf of ${\mathcal O}^{*}_V$ defined
as

\begin {itemize}
\item $J_q = ({\mathcal O}^{*}_V)_q$ if $q\notin C$.
\item $J_q = \{\phi \in ({\mathcal O}^{*}_V)_q; \phi|_C \simeq 1\}$ if $q\in C$.
\end{itemize}

We have then the short exact sequence
$$
1\rightarrow J \rightarrow {\mathcal O}^{*}_V \rightarrow
{\mathcal O}^{*}_V/J \rightarrow 1;
$$
we remark that ${\mathcal O}^{*}_V/J$ can be taken as
${\mathcal O}^{*}_C$.

In order to have the surjectivity stated above, we need
$H^2(V,J)=0$. Since the exponencial map gives an isomorphism between
$\mathcal I_C$ and $J$, it is enough to have $H^2(V, \mathcal
I_C)=0$. \qed

\section{\bf Constructing Foliations}

We are able now to prove the Theorem stated in the Introduction.

\medskip
Let $Y: C\rightarrow TS|_C$ be the line field constructed in Corollary 1.
Let $\{U_i\}$ be a covering of $C$ and  $\tilde U_i$ be an open set such that $\tilde U_i\cap C= U_i$. In each $\tilde U_i$ we choose a 1-form $\omega_i$ satisfying $ker (\omega_i(p))=Y(p)$ when $p\in U_i$. We may take coordinates  $(x_i,y_i) \in \tilde U_i$ as to have $U_i=\{y_i=0\}$ and $\omega_i=dy_i-x_i^{n_i}dx_i$(remember that the possibility  $n_i=0$ is allowed). We remark that $\omega_i|_{U_i\cap U_j} = f_{ij}\,\omega
_j|_{U_i\cap U_j}$ whenever $U_i\cap U_j\neq {\emptyset},
\,\,\,f_{ij}\in Z^1(\{U_i\},\mathcal O_{C}^{*}$). 
%A simple
%calculation gives $$ f_{ij}(x_j)= \frac
%{x_i^{n^i}}{x_j^{n^j}}\phi^{'}(x_j)
%$$
%\bigskip
We denote by $L=\{F_{ij}\}$ the line bundle over $V$ whose restriction
to $C$ is  defined by the transition functions $\{f_{ij}\}$ (Lemma 2); we have
$$
L|C= \mathcal O(D) \otimes KC^{*},
$$ where 
$D=
\sum_{i=1}^{l}n_ip_i$. The boundary $\delta \{\omega_i\}$ computed
in $Z^1(S,\Omega^1_S\otimes L)$ belongs effectively to
$Z^1(S,\mathcal I_C\cdot\Omega^1_S\otimes L)$, where
$\Omega^1_S$ is the sheaf of germs of holomorphic 1-forms of $S$.

\bigskip
We claim that $H^1(S,\mathcal I_C\cdot\Omega^1_S\otimes 
L)=0$. As discussed before, we need that $\forall \nu \geq 1$
$$
H^0(C,KC\otimes NC^{\nu}\otimes({\Omega^1_S \otimes 
L})^{*}|C)=0
$$

\noindent which depends on  $$ H^0(C,KC^2 \otimes NC^{\nu}\otimes
\mathcal O(-D)\otimes TC) =0 \,\,\,\forall \nu \geq 1
$$ \noindent and $$
H^0(C,KC^2 \otimes NC^{\nu}\otimes \mathcal O(-D)\otimes NC) =0
\,\,\,\forall \nu \geq 1;
$$  
\smallskip
\noindent both equalities are true since 
the Chern classes of the line bundles $ KC^2 \otimes NC^{\nu}\otimes
\mathcal O(-D)\otimes TC $ and $KC^2 \otimes NC^{\nu}\otimes
\mathcal O(-D)\otimes NC$ are negative due to the hypothesis.

\bigskip
It follows that there exists a 0-cocycle $\{\eta_i\}\in
H^0(\tilde U_i,\mathcal I_C\cdot\Omega^1_S\otimes L)=0$ such that

$$
\omega_i - F_{ij}\,\omega_j= \eta_i - F_{ij}\,\eta_j
$$

\bigskip

\noindent and the foliation we look for is defined by the 1-form $$
\{\omega_i - f_{ij}\,\eta_i\}\in H^0(V, \Omega^1_S\otimes  L).
$$  \qed

\begin{coro}Let $C\hookrightarrow S$ be an embedding of the curve $C$ into the surface $S$ such that $C\cdot C<0$. Then there exists a regular holomorphic foliation defined in a neighborhood of $C$.
\end{coro}

%%%%%%%%%%%%%%%%%%%%%%%%%%%%%%%%%%%%%%%%%%%%%%%%%%5
\section{\bf Examples}
%\begin {PROP}
%\label {atmostone} Assume $c(NC)>2-2g+\sum n_i$. There exists at
%most one  line field along $C$ with the tangency divisor $D=\sum
%n_ip_i$.
%\end{PROP}  \proof According to Proposition 2, it is enough to guarantee that the line
%bundle $\mathcal O(D)\otimes NC^*\otimes TC$ has no sections, which
%is certainly the case when $c(\mathcal O(D)\otimes NC^*\otimes
%TC)<0$. \qed
\begin{figure}[t]
\begin{center}
%\includegraphics[width=0.4\textwidth]{figures/weierstrass.pdf}
%\caption{Lamniscate: $(x^2+y^2)^2=x^2-y^2+t,\ t=-0.2,-0.1,0,0.1,0.3,1$}
\includegraphics{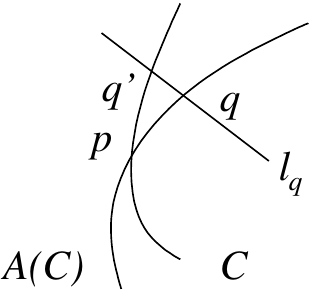}
\caption{}
\label{CAC}
\end{center}
\end{figure}
\begin{EX} 
\label{example1}\rm {\it A plane smooth projective curve $C$ different from the projective
line does not have a transverse holomorphic line field} (this is a
particular case of a theorem of Van de Ven (\cite {VV})). In fact,
suppose $Y$ is a transverse holomorphic line field defined along
$C$. We consider a holomorphic automorphism $A$ of the plane close
to the Identity which fixes some point $p\in C$ and such that
$(A_{*}Y)(p)\neq Y(p)$; the line field  $Y_A = A_{*}Y$ is of course
transverse to $A(C)$. Given $q\in A(C)$, we denote as $l_q$ the
projective line tangent to $Y_A(q)$ at $q$.    We may therefore
induce along $C$ a new holomorphic line field $Z\neq Y$ in the
following way: given $q\in A(C)$ take $q^{\prime}= l_q\cap C$ (the
intersection is taken in a small neighborhoood of $C$); then
$Z(q^{\prime})$ is  the tangent line to $l_q$ at the point
$q^{\prime}$, see Figure \ref{CAC}. Since $Z(p)= Y_A(p)\neq Y(p)$ and $Z$ is transverse to
$C$, we get a contradiction  with the Proposition \ref{atmostone}
(notice that $c(NC)=d^2 $ is greater than $3d- d^2 = 2-2g$ when
$d={\rm degree}(C)>1$).

A different, "foliated" argument goes as follows: suppose that there exists
a holomorphic line field transverse to $C$; this is a line bundle $L$ over $C$, as we have seen in Section 1 For each $p\in C$ the line $l(p)$ of $L$ is associated to
a projective line $l_p$ of the projective plane (this projective line passes through $p$ with direction given by $l(p)$). We take some
Riemannian metric in ${\mathbb P}^2$; since there are neighborhoods of $C$ in $L$ and of $C$ in ${\mathbb P}^2$ which are $C^{\infty}$ diffeomorphic (as line bundles) , for a small $\eta$ the discs
centered at the points of $C$ of radius $\eta$ and contained in the
projective lines $\{l_p\}_{p\in C}$ are disjoint, so they  form a holomorphic fibration. We
pick up  a non-constant meromorphic function in $C$ and extend it to
a neighborhood of $C$ as a constant along each fiber. This is a
meromorphic function that can be extend to all of ${\mathbb P}^2$
since the complement of $C$ is a Stein surface. We observe that the
extension is constant along each projective line $l_p$. The only
possibility is that these projective lines form a pencil   issued
from some point of the plane.

Proposition \ref{25mar2011} implies that for a curve with sufficiently positive self-intersection we have at most one transverse holomorphic line field. The above example shows that such a 
line field may not exist at all. Note that the construction of transverse line fields presented in Proposition \ref{25mar11} is done under the hypothesis that the curve has sufficiently negative
self-intersection.  
\end{EX}

\smallskip
\begin{EX} \rm The Proposition \ref{atmostone} is useful to get examples of non-existence of
certain regular foliations when the self-intersection of $C$ is not
sufficiently negative. In order to see this, let us consider a pair
$C\hookrightarrow S$ obtained by the following procedure:

\begin{enumerate}
\item we blow up the origin 0 of the polydisc $\Delta \subset \mathbb C^2$, introducing
an exceptional divisor; we choose the point in this divisor which
belongs to the strict transform of $\{y=0\}$ and blow up again. We
keep doing this in order to get a chain of projective lines
$E_1,...,E_{m-1}$ of self-intersection $-2$ and a last projetive
line $E_m$ of self-intersection $-1$; there is a holomorphic
projection $\pi$ from the resulting surface $\tilde \Delta$ to
$\Delta$, which collapses $E_1 \cup \dots \cup E_m$ to 0, and which
is an isomorphism from the complement of this divisor to
$\Delta\setminus \{(0,0)\}$. Denote by   $q\in E_m$ the point which
belongs to the strict transform of $\{y=0\}$ and take the
$u$-coordinate along $E_m$ in order  to have $\pi(x,u)=(x,ux^m)$. We take
also a polydisc $V= \{x,u); |x|<1, |u|< \epsilon\}$, for a small
$\epsilon$, around $(x,u)=(0,0)=q\in E_m$.
\item let us consider a line bundle over a compact, holomorphic, smooth curve $\tilde
C$ whose self-intersection satisfies $\tilde C \cdot \tilde C>
2-2g$; we select some point in $\tilde C$ and introduce coordinates
$(\tilde x,\tilde u)$ in a neighborhood $W$ of this point as to have
$\{\tilde x= {\rm const}\}$ contained in the linear fiber through
$(\tilde x,0)\in \tilde C$ for every $\tilde x$.
\item finally we glue $W$ to $V$ by means  of a holomorphic diffeomorphism $\Phi:W
\rightarrow V$ in order to get a holomorphic surface $\tilde S$
containing $E_1\cup \dots \cup E_m\cup \tilde C$ as a divisor whose
components have the self-intersection numbers described above;
$\Phi$ must send $(\tilde x,\tilde u)=(0,0)$ to $(x,u)=(0,0)=q$, the
$\tilde x$-axis into the $x$-axis and the $\tilde u$-axis
transversely to the $u$-axis. We remark that $\tilde C$ has a unique
field $\tilde {\mathcal L}$ of transversal lines because $\tilde C
\cdot \tilde C>2-2g$; by construction the line ${\tilde {\mathcal
L}}_q$ is different from $T_qE_m$.
\end{enumerate}

We blow down $E_1\cup \dots \cup E_m$ to $p=(0,0)\in \Delta$ and get
a surface $S$ with an embedded curve $C$ such that $C\cdot C
>m+2-2g$ and $p\in C$.

\medskip
\it {We claim that there exists no regular foliation $\mathcal F$ in
$S$ transverse to $C\setminus \{p\}$  with order of tangency  $0\leq
n \leq m-1$ at $p$ }.
\medskip
\rm {Otherwise after blowing up $a$ times as explained before
starting at $p$, we would get a foliation $\tilde {\mathcal F}$
transverse to $\tilde C$ and having $E_m$ as a leaf. Each leaf
$\tilde F_s$ through $s\in {\tilde C}$ has $\tilde L_s$ as tangent
line at $s\in \tilde C$; but this property is not verified at the
point $q\in \tilde C \cap E_m$}.

We remark that the particular case $m=1$ gives examples of
embeddings $C\hookrightarrow S$ such that $C\cdot C >3-2g$ without
transversal foliations to $C$; in particular, there is no
 neighborhood of $C$ in $S$ which is
(holomorphically) diffeomorphic to a neighborhood of $C$ in the total space of its
normal bundle.
\end{EX}
\bigskip

\section{\bf Plane curves and line fields}
We develop here Example \ref{example1} in order to understand the
role of tangencies. Let us consider in $\mathbb P^2$ a smooth
algebraic curve $C$ of degree $d$ and a holomorphic line field $X$
along $C$. We have then a holomorphic map ${\phi}_X:
C\longrightarrow \check{\mathbb P}^2$  
  defined as $\phi_X(p)= X(p)\in \check{\mathbb P}^2$; its image is
an algebraic curve $\check{X}\subset \check{\mathbb P}^2$. Let us denote
by $l\in {\mathbb N}$ the degree of $\phi$ as a map from $C$ onto $\check{X}$.

For instance, let us suppose that $X$ is induced by a pencil of
lines issued from some point $b\in \mathbb P^2$. Then $\check{X}$ is
a line in $\check{\mathbb P}^2$ and $\phi_X$ has degree $d$ or $d-1$
according to  $b\in C$ or $b \notin C$ (in this last case, $X(b)$ is
the tangent line to $C$ at $b\in C)$. We have then $tang(X,C)=d^2-d$
or $tang(X,C)=d^2-d-1$.

\begin{PROP} $tang(X,C)= l.deg(\check{X}) + d\,{^2} - 2d$.
\end{PROP}
\proof We consider $\mathbb P(T{\mathbb P}^2|_{C})$, which is a
$\mathbb P^1$-bundle over $C$ with the section $\mathbb P(TC)$. The
vector bundle $T{\mathbb P}^2|_{C}$ may be described by the
following transition maps:

$$ x_{\alpha} = \xi_{\alpha\beta}(z_{\beta})x_{\beta} + \eta
_{\alpha\beta}(z_{\beta})y_{\beta},\,\,\,\,\, y_{\alpha}=
c_{\alpha\beta}(z_{\beta})y_{\beta}
$$
where $(x_{\beta},y_{\beta})$ are coordinates for $T{\mathbb
P^2}|_{C}$ at the point of $C$ of coordinate $z_{\beta}$,
\,\,$z_{\alpha} = g_{\alpha\beta}(z_{\beta}),\,\,
\xi_{\alpha\beta}(z_{\beta})= g^{\prime}_{\alpha\beta}(z_{\alpha})$
and $\{c_{\alpha\beta}\}$ defines the normal bundle to $C$ in
$\mathbb P^2$.

In order to get the transition functions of $\mathbb P(T{\mathbb
P}^2|_{C})$, we put $u_{\beta}=x_{\beta}/y_{\beta}$ and $t_{\beta}=
y_{\beta}/x_{\beta}$; then

$$
u_{\alpha}=
\dfrac{\xi_{\alpha\beta}(z_{\beta})}{c_{\alpha\beta}(z_{\beta})}u_{\beta}
+ \dfrac {\eta
_{\alpha\beta}(z_{\beta})}{c_{\alpha\beta}(z_{\beta})}
$$
and $$ t_{\alpha}=
\dfrac{c_{\alpha\beta}(z_{\beta})t_{\beta}}{\xi_{\alpha\beta}(z_{\beta})
+ \eta _{\alpha\beta}(z_{\beta})t_{\beta}}
$$

Let us consider the line field $X$ as a section of $\mathbb
P(T{\mathbb P}^2|_{C})$; we choose also a generic pencil of lines
$P$. In the $u$-coodinates, we have

$$
X_{\alpha}=
\dfrac{\xi_{\alpha\beta}(z_{\beta})}{c_{\alpha\beta}(z_{\beta})}X_{\beta}
+ \dfrac {\eta
_{\alpha\beta}(z_{\beta})}{c_{\alpha\beta}(z_{\beta})}
$$
and $$ P_{\alpha}=
\dfrac{\xi_{\alpha\beta}(z_{\beta})}{c_{\alpha\beta}(z_{\beta})}P_{\beta}
+ \dfrac {\eta
_{\alpha\beta}(z_{\beta})}{c_{\alpha\beta}(z_{\beta})}
$$

\medskip
The intersection number of both sections $X$ and $P$ with $\mathbb
P(TC)$ will be denoted by $Poles(X)$ and $Poles(P)$; of course
$tang(X,C)= Poles(X)$ and $Poles(P)=d^2-d$.

\bigskip
>From the formulae above we see that $\{X_{\alpha}- P_{\alpha}\}$ is
a section of the linear bundle given by the cocycle $\{\dfrac{\xi_{\alpha\beta}(z_{\beta})}{c_{\alpha\beta}(z_{\beta})}\}$,
which is $TC\otimes NC^{*}$. Consequently:

$$
Zeroes(X-P)- Poles(X-P)= -2d^2+3d
$$

\medskip
Therefore $Poles(X)=Zeroes(X-P) - Poles(P)+ 2d^2 + 3d$. Now since
$Poles(P)=d^2-d$ and $Zeroes(X-P)=l.deg(\check{X})$, we get finally
$$ tang(X,C)=l.deg(\check{X})+ d^2-2d.
$$
\qed

\bigskip
\begin{coro}$tang(X,C)\geq (d-1)^2$
\end{coro}

This Corollary gives another explanation why a
  a smooth, plane algebraic curve $C$ of
degree greater than one has no transversal holomorphic line field;
consequently a neighborhood of $C$ can not be linearized. 

We see also  that {\it if we blow up at $d^2-2d$
different points of $C$, the resulting curve $\hat C$ has not a
linearizable neighborhood as well}. In fact, a tranversal
holomorphic line field to $\hat C$ corresponds to a holomorphic line
field along $C$ with at most $d^2-2d$ points of ordinary tangency,
which is not possible.

\end{document}